\documentclass[12pt]{amsart}
\usepackage{amsmath,amscd,amsthm,amsfonts, amssymb,amsxtra}
\usepackage[all]{xy}
\SelectTips{cm}{}
\newtheorem{thm}{Theorem}[section]  
\newtheorem*{un-no-thm}{Theorem}
\newtheorem{cor}[thm]{Corollary}     % Numbered along with thm
\newtheorem{lem}[thm]{Lemma}         % Numbered along with thm
\newtheorem{prop}[thm]{Proposition}  
\newtheorem{bigthm}{Theorem}
\theoremstyle{definition}
\theoremstyle{definition}
\theoremstyle{definition}

\theoremstyle{remark}
\newtheorem{rem}[thm]{Remark}

\newtheorem*{acks}{Acknowledgements}
 
\newtheorem*{out}{Outline}

\newtheorem*{notation}{Notation}

\begin{document}
\title{Poincar\'e complex diagonals}
\date{\today}
\author{John R. Klein}
\address{Wayne State University, Detroit, MI 48202}
\email{klein@math.wayne.edu}
\begin{abstract} Let $M$ be a
Poincar\'e duality space of dimension $d\ge 4$.
In this paper we describe a complete obstruction to 
realizing the diagonal map $M \to M {\times} M$
by a Poincar\'e embedding. The obstruction group
depends only on the fundamental group and the parity
of $d$.
\end{abstract}
\thanks{The author was partially supported by the NSF}
\maketitle%\pagestyle{fancy}
\setlength{\parindent}{15pt}
\setlength{\parskip}{1pt plus 0pt minus 1pt}
\def\Top{\bold T\bold o \bold p}
\def\wTop{\text{\rm w}\bold T}
\def\wT{\text{\rm w}\bold T}
\def\Sp{\bold S\bold p}
\def\vo{\varOmega}
\def\vs{\varSigma}
\def\smsh{\wedge}
\def\flush{\flushpar}
\def\id{\text{id}}
\def\dbslash{/\!\! /}
\def\codim{\text{\rm codim\,}}
\def\:{\colon}
\def\holim{\text{\rm holim\,}}
\def\hocolim{\text{\rm hocolim\,}}
\def\hodim{\text{\rm hodim\,}}
\def\hocodim{\text{hocodim\,}}
\def\Bbb{\mathbb}
\def\bold{\mathbf}
\def\Aut{\text{\rm Aut}}
\def\cal{\mathcal}
\def\frak{\mathfrak}
\def\Sec{\text{\rm sec}}
\def\Secst{\text{\rm sec}^{\text{\rm st}}}
\def\map{\text{\rm map}}
\def\stableto{\mapstochar \!\!\to}

\setcounter{tocdepth}{1}
%\tableofcontents
%\addcontentsline{file}{sec_unit}{entry}
%\endtableofcontents

\section{Introduction \label{intro}}

A distinguishing feature of smooth manifolds 
is that they come equipped with a preferred
bundle of tangent vectors. If $M$ is a smooth manifold
with tangent bundle $\tau_M$, then the diagonal map
\begin{align*}
M  &\to M \times M \\
x &\mapsto (x,x)
\end{align*}
is a smooth embedding possessing a tubular neighborhood 
diffeomorphic to the total space of $\tau_M$. Moreover, the
set of tubular neighborhoods of the diagonal can be topologized
in such a way that the resulting space is contractible.

In this note we shall consider consider a homotopy theoretic analog
of the above in which $M$ is now replaced by a
space satisfying Poincar\'e
duality. There are two parts 
question to consider: existence and uniqueness.
Here we will be concerned with the existence problem.

Forty years ago, Spivak \cite{Spivak} proved that a Poincar\'e 
space possesses a ``stable tangent bundle'' in the sense of stable spherical
fibration theory. By Atiyah duality \cite{Atiyah}, 
the Spivak tangent fibration is
{\it stably} unique up to contractible choice. 
Furthermore, elementary obstruction theory
shows that a stable spherical fibration over a 
$d$-dimensional Poincar\'e space
lifts to an unstable $(d{-}1)$-spherical fibration. 
However, it isn't unique:
the number of such lifts is countably infinite when $d$ is even and
precisely two when $d$ is odd.

Given a choice of unstable lift of
the Spivak fibration to a $(d{-}1)$-spherical fibration $\tau\:S(\tau) \to M$,
one can then ask whether
it appears as normal data for an embedding of the diagonal, where
now by ``embedding'' we mean a {\it Poincar\'e embedding} in the
sense the diagonal map extends to a homotopy equivalence
$$
M \times M \,\, \simeq \,\, D(\tau) \cup_{S(\tau)} C \, ,
$$
in which $D(\tau)$ is the mapping cylinder of $\tau$ and $(C,S(\tau))$
is a Poincar\'e pair. It may be the case that for some choice of $\tau$
there is such a Poincar\'e embedding. Our main result will be to provide
a complete obstruction in dimensions at least four. 
It will turn out that the
obstruction group depends only on the fundamental group of $M$
and the parity of $d$.

To describe the obstruction group,
let $\pi$ be the fundamental group of $M$ and  let
$\bar\pi$ be its set of conjugacy classes. 
We denote the conjugacy class of
a group element $x$ by $\bar x$. Let ${\Bbb Z}[\bar \pi]$
be the free abelian group with basis $\bar \pi$ .
Define an involution
on ${\Bbb Z}[\bar \pi]$ by
$$
\bar x \mapsto {(-1)}^d \bar x^{-1}, \qquad x \in \pi
$$
and extending linearly. Let 
$$
Q_d(\pi) \,\, := \,\, {\Bbb Z}[\bar \pi]_{\Bbb Z_2}
$$
be the resulting abelian group of coinvariants.
When $\pi = e$ is the trivial group, $Q_d := Q_d(e)$ is infinite
cyclic or cyclic of order two depending on whether $d$ is even or odd.
Let 
$$
\tilde Q_d (\pi) \,\, = \,\, \text{coker}(Q_d \to Q_d(\pi))
$$
be the cokernel of the evident homomorphism $Q_d \to Q_d(\pi)$.

\begin{bigthm} \label{first} Associated with
each connected Poincare duality space $M$ of 
dimension $d$ and fundamental group $\pi$, there is an obstruction
$$
\mu_M \in \tilde Q_d(\pi)\, ,
$$
which depends only on the homotopy
type of $M$ and which vanishes when the diagonal $M\to M\times M$ 
admits a Poincar\'e embedding.

If $d \ge 4$ and $\mu_M$ is trivial, then 
there is such a Poincar\'e embedding.
\end{bigthm}

{\flushleft \it Remarks:}

{\flushleft (1).} The group
$\tilde Q_d(\pi)$ is trivial if and only if $\pi$ is. 
Consequently, if $M$ is simply connected, 
there is no obstruction and $M$ diagonally 
Poincar\'e embeds, recovering \cite[Cor.\ H]{Klein_compression}.
\medskip

{\flushleft (2).} Aside from manifolds, there 
are some special cases with non-trivial $\pi$ which give trivial 
$\mu_M$. For example, according to Byun \cite{Byun},
when $M$ is obtained by gluing two compact 
smooth manifolds along their boundaries 
using a homotopy equivalence (i.e, a two-stage patch space) 
and  the square-root closed condition is satisfied by the 
fundamental group of the boundary, then $M$ admits
a diagonal embedding.
\medskip

{\flushleft (3).} I do not know if $\mu_M$ is 
always trivial, but I suspect it needn't be. 
If we view uniqueness as a kind of relative
existence, then some evidence against the global triviality
of $\mu_M$ is provided by
recent work of Longoni and Salvatore \cite{LS} which essentially
exhibits  two distinct concordance classes of 
diagonal Poincar\'e embeddings
of the lens space $L(7,1)$.
\medskip

{\flushleft (4).} Dupont \cite{Dupont1}, \cite{Dupont2},
Sutherland \cite{Sutherland} and Byun \cite{Byun2}
showed that the Spivak tangent fibration always admits a preferred
 $(d-1)$-spherical lift. If $\mu_M = 0$, one also obtains
a preferred lift. The relationship between the two lifts 
is still unclear to me. 
\medskip

To find the unquotiented group $Q_d(\pi)$, we were
heavily influenced by the paper of Hatcher and Quinn
\cite{Hatcher-Quinn}, which
gives a bordism theoretic approach to deciding, in 
the metastable range,
when an immersion of smooth manifolds is regularly homotopic
to an embedding. Their obstruction
was defined geometrically by placing the immersion
in self-transverse position, taking the resulting
double point manifold and considering it as a bordism class.
In the special case of the diagonal map
of  a smooth manifold $M$, the Hatcher-Quinn bordism group
is isomorphic to $Q_d(\pi)$ (of course, since the diagonal 
of a smooth manifold is an embedding,
its Hatcher-Quinn invariant is trivial). 

In the Poincar\'e category
we cannot use the Hatcher-Quinn approach because
the requisite transversality theory is lacking. 
Nevertheless, there is
a way to work {\it externally} to define a homotopy
theoretic version of the Hatcher-Quinn invariant. A crude version of 
the idea is this: think of an immersion $P \to N$ as a ``cycle.'' 
It is then Poincar\'e dual to a ``cocycle'' 
which in this case is a {\it stable} 
Thom-Pontryagin map
$N \to P^\nu$, where $P^\nu$ is the Thom space of the stable normal
bundle (= stable spherical fibration in the Poincar\'e setting). 
If the immersion were an embedding, this map would desuspend to 
an actual map. This suggests
that the obstruction we seek should be the one which realizes when
the map desuspends. It is well-known that the Hopf invariant 
provides such an obstruction, once an 
unstable representative of $\nu$ has been chosen. 
In this paper, we will use
a $\pi$-equivariant version of the Hopf invariant 
currently being developed by Crabb and Ranicki 
to capture some additional information lost
by the ordinary Hopf invariant. 

In addition to the above, the other problem
one needs to consider is the indeterminacy 
associated with choosing different
lifts of $\nu$ to an unstable $(n{-}p{-}1)$-spherical fibration.
In the case of the diagonal map, the indeterminacy vanishes when
one takes the quotient of $Q_d(\pi)$ by the subgroup $Q_d$. 
This leads to the definition of $\mu_M$.

\begin{out} \S2 consists of preliminary material.
In \S3 we define variants of the Thom-Pontryagin collapse
of a Poincar\'e embedding. In \S4 we consider the problem
of compressing Poincar\'e embeddings into one dimension
beyond the stable range. \S5
develops the equivariant stable Hopf invariant and uses it
to reinterpret the results of \S4. In \S6 we consider
of Poincar\'e embeddings of the diagonal for fixed choice of 
$(d{-}1)$-spherical lift of the stable 
tangent fibration.  \S7  explains how the 
set of equivalence classes of $(d{-}1)$-spherical
lifts of the stable tangent fibration form a $Q_d$-torsor. 
In \S8 we prove Theorem \ref{first}.
\end{out}

\begin{acks} I first learned about the equivariant Hopf invariant
in conversations with Bill Richter. I am indebted to Andrew Ranicki
for explaining to me his work with Michael Crabb.
\end{acks}

\section{Preliminaries}

\subsection*{Spaces} The spaces of this paper are to be topologized
with the compactly generated topology.
A map $X\to Y$ of based spaces is a {\it weak equivalence} if it is 
a weak homotopy equivalence. A weak equivalence is written with 
the symbol `$\overset \sim \to$', and we use the notation $X\simeq Y$
to indicate that $X$ and $Y$ are equated by a chain of weak equivalences.
We write $[X,Y]$ for  based homotopy
classes of maps from $X$ to $Y$. 

A non-empty (possibly unbased) space $Z$ 
is {\it $0$-connected} if it is path connected.
It is $r$-connected with  $r > 0$, if it is $0$-connected and
$[S^k,Z]$ has one element for $k\le r$.
A map of non-empty (possibly unbased) spaces is
{\it $r$-connected} if its homotopy fiber, taken with respect to
all choices of basepoint, is $(r{-}1)$-connected.

A {\it stable map} from $X$ to $Y$ is a map $X\to QY$ where
$Q = \Omega^\infty \Sigma^\infty$ is the stable homotopy functor.
Alternatively, it is colimit under suspension of maps of the form
$\Sigma^j X \to \Sigma^j Y$. Stable maps can be composed.
We write $\{X,Y\}$ for the  abelian group of stable 
homotopy classes of maps from $X$ to $Y$. A stable
map is denoted with `$\,\,\stableto$.' We sometimes write $X \simeq_s Y$
to indicate that $X$ and $Y$ are related by a chain of
stable weak equivalences.

If $G$ is a group and $X$ and $Y$ are now based $G$-spaces, then
an equivariant map $X\to Y$ is a weak equivalence if it is a weak homotopy
equivalence of underlying spaces. We write $[X,Y]_G$ to indicate 
equivariant homotopy
classes of equivariant based maps, and  $\{X,Y\}_G$ for the
equivariant stable homotopy classes. We note that the kind of equivariant
homotopy we are considering here is the naive one in which the group is 
acting trivially on suspension coordinates, and all objects are
weak equivalent to ones with free actions.

\subsection*{Poincar\'e spaces} 
A finitely dominated space $X$ is a {\it Poincar\'e duality space} 
of (formal) dimension $d$ if there exist
a bundle of coefficients $\mathcal L$ which is locally isomorphic
to $\Bbb Z$, and a class $[X] \in H_d(X;{\mathcal L})$
such that the associated cap product homomorphism
$$
\cap [X]\:H^*(X;M) \to H_{d{-}*}(X;{\mathcal L}\otimes M)
$$
is an isomorphism in all degrees. Here, $M$ denotes any 
bundle of coefficients (cf.\  \cite{Wall_PD}, \cite{Klein_PD}).
The pair
$({\cal L},[X])$
is uniquely defined up to unique isomorphism; ${\cal L}$ is called
the {\it orientation bundle} and $[X]$ a {\it fundamental class}.

Similarly, one has the definition of Poincar\'e space
$X$ {\it with boundary} ${\partial X}$ (also called a
{\it Poincar\'e pair} $(X,\partial X)$).
Here, one assumes both $X$ and ${\partial X}$ are finitely
dominated and there is a fundamental class 
$[X] \in H_d(X,\partial X;{\mathcal L})$
such that 
$$
\cap [X]\:H^*(X;M) \to H_{d{-}*}(X;\partial X; {\mathcal L}\otimes M)
$$
is an isomorphism.  Additionally, if $[\partial X]$ is
 the image of $[X]$ under the boundary homomorphism
$H_d(X,\partial X;{\cal L}) \to 
H_{d-1}(\partial X;{\cal L}_{|\partial X})$, one 
also requires 
$$
({\cal L}_{|\partial X},[\partial X])
$$ 
to equip $\partial X$ with the structure of a Poincar\'e space. 

Unless otherwise indicated, we usually work with 
{\it connected} Poincar\'e spaces. We say that a Poincar\'e space $X$ is
{\it closed} if it is without boundary.

\subsection*{Poincar\'e embeddings}
Let 
$$
f\: P \to N
$$ 
be a map of connected Poincar\'e duality spaces, where
\begin{itemize}
\item $P$ is closed, but $N$ need not be;
\item $P$ has dimension $p$, $N$ has dimension $n$; 
\item $p \le n -1$.
\end{itemize}
A {\it Poincar\'e
embedding} of a $f$ consists of an $(n{-}p{-}1)$-spherical 
fibration 
$$
\xi \:S(\xi) \to P
$$
and commutative diagram
$$
\xymatrix{
S(\xi) \ar[r]\ar[d] & C \ar[d] &\partial N \ar[l]\\
D(\xi) \ar[r]_{f_\xi}& N\, ,
}
$$
such that
\begin{itemize}
\item $D(\xi)$ is the mapping cylinder $\xi$ (in particular,
$D(\xi) \simeq P$), and $f_\xi$ is the composite
$$
D(\xi) \simeq P \overset f \to N\, .
$$
\item The square in the diagram is a homotopy pushout.
\item The composite $\partial N \to C \to N$ is the inclusion.
\item The image of a fundamental class $[N]$
under the composite
$$
H_n(N,\partial N;{\cal L}) \to H_n(N,C;{\cal L}) \cong 
H_n(D(\xi),S(\xi); {\cal L}_{|D(\xi)})
$$
gives $(D(\xi),S(\xi))$ the structure of a Poincar\'e space.

Similarly, the image of a fundamental class under
$$
H_n(N,\partial N;{\cal L}) \to H_n(N,D(\xi) \amalg \partial N;{\cal L}) 
\cong H_n(C,S(\xi) \amalg \partial N; {\cal L}_{|C})
$$
equips $(C,D(\xi)\amalg \partial N)$ with the structure of
a Poincar\'e space.
\end{itemize}

We warn the reader that to avoid notational clutter, we
are not necessarily assuming in the
above that $S(\xi) \amalg \partial N \to C$ is an
inclusion, and we are slightly abusing notation when writing
expressions such as $H_n(N,C;{\cal L})$
(the reader should substitute the appropriate
mapping cylinder in such cases).

The spherical fibration $\xi$ appearing above is called the {\it normal datum}
of the Poincar\'e embedding.
We say in this case that that we have
a {\it codimension zero} Poincar\'e embedding
of the map $f_\xi\: D(\xi)\to N$.

We will define concordance only in the codimension zero setting.
Note that a codimension zero Poincar\'e embedding is really just
a certain kind of object 
$$
C\in \Top_{S(\xi)\amalg \partial \to N}\, ,
$$
where the latter is the category of spaces which factorize
the given  map $S(\xi)\amalg \partial \to N$. A morphism
in this category is just a map of spaces which is compatible
with the given factorizations. We say that such a morphism 
is a weak equivalence if it is a weak homotopy equivalence
of underlying topological spaces. Given two Poincar\'e embeddings
$C,C'$, considered as objects in the above category, we say that
they are {\it concordant} if there is a chain of weak equivalences
connecting them.

\subsection*{Poincar\'e Immersions}

Let $f\: P^p \to N^n$ be as above, and suppose we are given
an $(n{-}p{-}1)$-spherical fibration $\xi\:S(\xi) \to P$.
By analogy with Smale-Hirsch theory, a {\it Poincar\'e immersion} of 
$(f,\xi)$ is
a choice of (codimension zero) Poincar\'e embedding of the map
$$
\begin{CD}
f_{\xi\oplus \epsilon^j} \: D(\xi\oplus \epsilon^j) = 
D(\xi)\times D^j @> f_\xi \times \text{id}_{D^j} >> N \times D^j
\end{CD}
$$
for some $j \ge 0$. Here $\epsilon^j$ denotes the trivial
$(j{-}1)$-spherical fibration over $P$ and $\oplus$ is the fiberwise
join operation.

In \cite{Klein_immersion}, it was shown that $(f,\xi)$ immerses if and only
if there is a fiber homotopy equivalence
$$
\nu_P \oplus \xi \,\,  \simeq \,\, f^*\nu_N \, ,
$$
where $\nu_P$ and
$\nu_N$ are the Spivak fibrations of $P$ and $N$ respectively.

\section{The collapse}

\subsection*{The classical collapse}
Given $(f,\xi)$ as above, and
a Poincar\'e embedding
$$
\xymatrix{
S(\xi) \ar[r]\ar[d] & C \ar[d] &\partial N \ar[l]\\
D(\xi) \ar[r]_{f_\xi}& N\, ,
}
$$
we get maps of pairs
$$
\begin{CD}
(N,\partial N) @>>> (N,C) @<<< (D(\xi),S(\xi))
\end{CD}
$$
which yields {\it collapse map}
$$
\begin{CD}
N/\partial N @>>> N/C @< \sim << P^\xi
\end{CD}
$$  
upon taking quotients, where $P^\xi = D(\xi)/S(\xi)$ is
the {\it Thom space} of $\xi$.

This gives a well defined homotopy class 
$$
c_f \in [N/\partial N,P^\xi] \, .
$$

\subsection*{The equivariant version}
Let
$$
\tilde N \to N
$$
be a universal cover with $\pi$ its group of deck transformations.
Then $\pi$ is identified with the fundamental group of $N$ and
it acts freely  on $\tilde N$.  The orbit space 
$\tilde N/\pi$  coincides with $N$.

For any map $X \to N$, we define
$$
\tilde X \,\, = \,\, X\times^N \tilde N
$$
to be the fiber product of $X$ and $\tilde N$ along $N$. 
Then $\tilde X$ inherits a free $\pi$-action.

Returning again to our Poincar\'e embedding diagram, we pull all terms
back along $\tilde N$ to get a diagram of free $\pi$-spaces
$$
\xymatrix{
\tilde S(\xi) \ar[r]\ar[d] & \tilde C \ar[d] &\partial \tilde N \ar[l]\\
\tilde D(\xi) \ar[r]_{\tilde f_\xi}& \tilde N\, ,
}
$$
where $\partial \tilde N$ is shorthand for $\partial N \times^N \tilde N$.
Note that the square in this diagram is a homotopy pushout because
it comes about by pulling back a homotopy pushout along $\tilde N$.

Arguing just as in the classical case, we get a collapse map which
is now $\pi$-equivariant. So we obtain an equivariant homotopy class
$$
\tilde c_{f} \in [\tilde N/\partial \tilde N,\tilde P^\xi]_\pi \, ,
$$
where $\tilde P^\xi := \tilde D(\xi)/S(\xi)$
is the {\it equivariant Thom space}. 
We call this construction the
{\it equivariant collapse} of the Poincar\'e embedding.

If $\partial N = \emptyset$, then observe $\tilde N/\partial \tilde N$
is $\tilde N^+ = $ the union of $\tilde N$ with a disjoint basepoint.

\subsection*{The  equivariant stable collapse}
Given $f\: P \to N$, we consider the virtual spherical fibration
$$
f^* \tau_N - \tau_P
$$
over $P$. This formal difference is to be interpreted in
the Grothendieck group of spherical fibrations over $P$.

If $j$ is sufficiently large, the sum 
$\xi :=  \tau_N - \tau_P \oplus \epsilon^j$ will admit
a spherical fibration representative. It follows that
when $j$ is large, the composite
$$
\begin{CD}
P  @> f >> N \subset N \times D^j 
\end{CD}
$$
will admit a Poincar\'e embedding with normal 
datum $\xi$.

Then we get an equivariant collapse, which 
is an element of
$$
[\Sigma^j \tilde N/\partial \tilde N,P^\xi]_\pi\, .
$$
Letting $j$ tend to infinity, we obtain an equivariant stable homotopy
class
$$
\tilde c_f \in \{\tilde N/\partial \tilde N,\tilde P^{f^*\tau_N - \tau_P}\}_\pi\, ,
$$
where $\tilde P^{f^*\tau_N - \tau_P}$
is the {\it equivariant Thom spectrum} of the virtual spherical fibration
$f^*\tau_N - \tau_P$.

We call $\tilde c_f$ the
{\it equivariant stable collapse}; it only depends on the homotopy
class of $f\: P \to N$.

\section{Compression}

If $Y \to N$ is a map, i.e., $Y$ is a ``space over'' $N$, then
we can form the {\it unreduced fiberwise suspension}
$$
S_N Y\,\,  :=  \,\, N \times 0 \cup_{Y\times 0} Y\times [0,1] \cup_{Y\times 1} N\times 1 
$$
which is a space over $N\times D^1 \cong N\times [0,1]$. 
Let
$$
S^j_N Y 
$$
be its $j$-fold iterate; it is a space over $N\times D^j$.

Given a Poincar\'e embedding
$$
\xymatrix{
S(\xi) \ar[r]\ar[d] & C \ar[d] &\partial N \ar[l]\\
D(\xi) \ar[r]_{f_\xi}& N\, ,
}
$$
we apply fiberwise
suspension $j$-times 
to obtain a new Poincar\'e embedding
$$
\xymatrix{
S(\xi\oplus \epsilon^j) 
\ar[r] \ar[d] &S^j_N C \ar[d] & S^j_N \partial N \ar[l] \\
D(\xi \oplus \epsilon^j) \ar[r]_{f_{\xi\oplus \epsilon^j}}& N\times D^j\, .
}
$$
(Note that $S^j_N\partial N$ is $\partial (N\times D^j)$.)
This is the
{\it $j$-fold decompression} of the given Poincar\'e embedding.
\medskip

We wish to find obstructions to 
reversing this procedure up to concordance. That is,
we wish to know when there is a Poincar\'e embedding of $f\: P \to N$
whose which decompresses to this Poincar\'e immersion up to 
concordance.

To answer this question, we take the equivariant stable collapse
$$
\tilde c_f \in \{\tilde N/\partial \tilde N,
\tilde P^{f^*\tau_N - \tau_P}\}_\pi
$$
and ask whether it comes from an {\it unstable class}.

More precisely, assume that $2p = n$. Then by 
classical obstruction theory,
there exists an unstable $(p{-}1)$-spherical fibration $\xi$ over $P$
which represents $f^*\tau_N - \tau_P$. Note that
$\xi$ isn't necessarily unique.

\begin{thm} \label{compression_result} 
Assume $n=2p$ and $p \ge 4$.  Choose a $(p{-}1)$-spherical
fibration $\xi$ representing $f^*\tau_N - \tau_P$.
Then $f$ admits a Poincar\'e embedding 
with normal datum $\xi$ if and only if $\tilde c_f$
destabilizes to an element of $[\tilde N/\partial \tilde N,
\tilde P^\xi]_\pi$.
\end{thm}

\begin{proof} This is a corollary
of \cite[th.\ B]{Klein_compression} and the proof
is essentially the same as that of \cite[th.\ E]{Klein_compression}
(for the proof, see \S4 of that paper).
\end{proof}

\section{The equivariant stable Hopf invariant}

Let $X$ and $Y$ be based spaces.
The stable Hopf invariant is a certain
function
$$
H\:\{X,Y\} \to \{X,D_2Y\}\, ,
$$
where $D_2 Y$ is the quadratic construction on $Y$.
This function is a natural transformation in $X$, and
hence, by the Yoneda lemma, it is determined by a map 
$$
QY \to QD_2 Y \, .
$$
For a stable map $f\: X\stableto Y$, the class
$H(f)$ is an obstruction to desuspending $f$
to an unstable homotopy class. It is a complete obstruction in the
metastable range.

In this paper, we will need to generalize 
the stable Hopf invariant to the $G$-equivariant setting for
a discrete group $G$. Let $X$ and $Y$ now be based $G$-spaces.
Then we will describe
a function
$$
\tilde H\: \{X,Y\}_G \to \{X,D_2Y\}_G \, ,
$$
which gives an obstruction to equivariantly 
destabilizing a stable $G$-map.

The construction of $\tilde H$ 
is basically the same as one of the well-known constructions of 
$H$; the details will  appear in a
forthcoming book of Crabb and Ranicki \cite{Crabb-Ranicki}
(see \cite{Crabb-Ranicki2} for a preliminary version). 
We will sketch their construction,
which uses G.\ Segal's  $\Bbb Z_2$-equivariant stable homotopy functor.

Let ${\cal U}$ be a complete universe of $\Bbb Z_2$-representations.
For example, we can take ${\cal U}$ to be a countable direct sum of
the regular representation. For $V \in {\cal U}$ we let $S^V$ denote its
one point compactification. Given $  Y$ as above, we set 
$$
Q_{{\Bbb Z}_2}(  Y\smsh   Y)\,\,  := \,\, 
\underset{V\in {\cal U}}{\text{\rm colim }} 
\text{map}(S^V,S^V\smsh   Y\smsh   Y)\, ,
$$
where the function space appearing above is given the
${\Bbb Z_2}$-action arising from conjugating maps.
Then $Q_{{\Bbb Z}_2}(  Y\smsh   Y)$ is a $(G\times \Bbb Z_2)$-space. 
Let
$$
Q_{\Bbb Z_2}(Y\smsh Y)^{\Bbb Z_2}
$$
denote the fixed point set of ${\Bbb Z}_2$ acting on $Q_{\Bbb Z_2}(Y\smsh Y)$.
This is a $G$-space.

If $Z$ is a based $G$-space, 
we let 
$$
\{Z,  Y\smsh   Y\}_G^{\Bbb Z_2}
$$
denote the abelian group of $G$-equivariant homotopy classes of $G$-maps
$Z \to Q_{\Bbb Z_2}(  Y\smsh   Y)^{\Bbb Z_2}$.
Note that such a homotopy class
is determined by specifying a $V\in {\cal U}$ and a 
$(G\times {\Bbb Z}_2)$-equivariant map
$$
S^V \smsh  Z \to S^V \smsh  Y\smsh  Y\, 
$$ 
(however, not all homotopy classes need arise in this way).
In particular, the reduced diagonal 
$\Delta_{ Y}\: S^0 \smsh  Y \to  S^0\smsh  Y\smsh  Y$
is such a map.

Now, given a stable $G$-map 
$$
f\:  X \stableto  Y \,
$$
we consider the difference
$$
\delta(f) := (f\smsh f)\Delta_X - \Delta_Y f \in 
\{ X, Y\smsh  Y\}_G^{\Bbb Z_2} \, .
$$

According to a formula of tom Dieck 
\cite[th.\ 2]{tomDieck} (cf.\ \cite[cor.\ A.3]{Crabb},  
\cite{Carlsson}, \cite[p.\ 203]{May_et_al}), 
there is a 
homotopically split fibration sequence of 
$G$-spaces
$$
QD_2  Y \overset{i}\to 
Q_{\Bbb Z_2}( Y\smsh  Y)^{\Bbb Z_2} 
\overset{j}\to Q 
Y\, .
$$
In particular, we have a (split) short exact sequence
$$
0\to \{ X,D_2  Y\}_G \overset{i_*}\to
\{ X, Y\smsh  Y\}_G^{\Bbb Z_2} \overset{j_*}\to
\{ X, Y\}_G\to 0 \, .
$$
One observes that the composite 
$$
j_*\delta(f)
$$
is trivial. Consequently, there is a unique element
$$
\tilde H(f) \in \{ X,D_2  Y\}_G
$$
such that $i_*\tilde H(g) =\delta(g)$\, . We call $\tilde H(f)$
the {\it $\pi$-equivariant stable Hopf invariant} of $f$.\footnote{Crabb
and Ranicki use the terminology {\it geometric Hopf invariant for 
$\tilde H(f)$}.}

\begin{notation} We will usually use 
$\tilde H$ when the group $G$ is understood. When the group
is ambiguous we resort to the notation $H_G$ to refer to
$\tilde H$. In particular, for the trivial group $e$, $H_e = H$
is the classical stable Hopf invariant.
\end{notation}

\subsection*{Desuspension}
Assume $X$ has the equivariant weak homotopy type of a $G$-space
obtained from a point by attaching free $G$-cells of dimension
at most $k$. Assume 
$Y$ is $r$-connected.

\begin{prop} Assume in addition $k \le 3r+1$. Then $\tilde H(f) = 0$
if and only if $f\in  \{X,Y\}_G$ desuspends to an  
element of $[X,Y]_G$.
\end{prop}

\begin{proof} (Sketch). First consider The case when 
$X = X_0\smsh (G)_+$, where $X_0$ is an unequivariant
based space. Then a stable $\pi$-map 
$f\:X\stableto Y$ is equivalent to specifying a stable map 
$f_0\:X_0\stableto Y$, and $\tilde H(f)$ is identified with
$H(f_0)$. The result now follows in this case from the fact
that, in the range $k \le 3r+1$, 
$H(f_0) = 0$ if and only if $f_0$ is represented by an unstable
map (cf.\ \cite{Milgram}). 

In the general case, one argues inductively by skeleta in
the equivariant CW decomposition for $X$. The previous paragraph 
amounts, more-or-less, to the inductive step.
\end{proof}

\begin{cor} \label{EHPsequence} With respect to the above assumptions,
the sequence of based sets
$$
\begin{CD}
[X,Y]_G @>>> \{X,Y\}_G @>\tilde H >> \{X,D_2 Y\}_G
\end{CD}
$$
is exact.
\end{cor}

\subsection*{The composition formula}
We will need to know how $\tilde H$ behaves with
respect to compositions.
Suppose $f\: X\stableto Y$ and $g\: Y\stableto Z$ are stable
$G$-maps. 

\begin{prop} \label{composition}
$$
\tilde H(g\circ f) = \tilde H(g)\circ f + D_2(g)\circ f \, .
$$
\end{prop}

This will be proved in \cite{Crabb-Ranicki} (cf.\
\cite{Crabb-Ranicki2}).
The formula is well-known when $G$ is the trivial group. When
$G$ is non-trivial one basically copies the proof for trivial $G$
noticing that all constructions are equivariant.

\subsection*{Extension of groups}
Let $K\to G$ be a group homomorphism, and suppose
$f\: X\stableto Y$ is an equivariant stable map of based 
$K$-spaces.
We assume without loss in generality that $X$ and $Y$ 
are built up from a
point by attaching free $K$-cells.
Then we have a  stable map of $G$-spaces
$$
\begin{CD}
f_{K\uparrow G} = f\smsh \text{id}_{G_+}\: X\smsh_K (G_+) \, 
\stableto Y\smsh_K (G_+) 
\end{CD}
$$
given by {\it inducing $f$ along $K\to G$.}
The inclusion $Y \to Y\smsh_K (G_+)$ is
$K$-equivariant and therefore induces a $K$-equivariant map
of quadratic constructions
$$
D_2Y \to  D_2 (Y\smsh_K (G_+))
$$
which admits a preferred extension to a $G$-map
$$
i_2 \:D_2Y\smsh_K (G_+) \to  D_2 (Y\smsh_K (G_+))  \, .
$$
The following lemma relates  
$H_G (f_{K\uparrow G})$ to $H_K(f)$.
Its proof, which
we will omit, is a straightforward
diagram chase.

\begin{lem}\label{induction} 
$$
H_G (f_{K\uparrow G}) \,\, = \,\, i_2 \circ (H_K (f) \smsh_K \text{\rm id}_{G_+}) \, .
$$
\end{lem}

\subsection*{Application to compression}
We now apply \ref{EHPsequence} to the problem of 
compressing Poincar\'e embeddings. As in the previous 
section, we are given a map $f\: P^p \to N^n$, with 
$P$ and $N$ Poincar\'e spaces, with $P$ closed and $n=2p$. Let $\pi$
be the fundamental group of $N$.

Then we have the associated equivariant stable collapse
$$
\tilde c_f\in \{\tilde N/\partial \tilde N, \tilde P^{f^*\tau_N -\tau_P}\}_\pi \, .
$$
Given a choice of unstable lift $\xi$ of $f^*\tau_N -\tau_P$
to an $p{-}1)$-spherical fibration,  we obtain an
explicit identification
\begin{equation}\label{thom_space}
\tilde P^{f^*\tau_N -\tau_P} := \tilde P^\xi \, .
\end{equation}
which equates the equivariant Thom spectrum on the left with
a representative equivariant Thom space on the right.

Using this identification, the equivariant stable Hopf invariant 
is realized as a function
$$
\tilde H^\xi\: \{\tilde N/\partial \tilde N,\tilde P^{f^*\tau_N -\tau_P}\}_\pi \, .
\to \{\tilde N/\partial \tilde N,D_2 \tilde P^{f^*\tau_N -\tau_P}\}_\pi \, ,
$$
in which are now indicating $\xi$ as a superscript the notation 
since it uses \eqref{thom_space}.

By \ref{compression_result} together with \ref{EHPsequence}, 
we  immediately obtain 

\begin{prop} \label{better_compression_result} 
Assume $n=2p$ and $p \ge 4$.
Then $f_\xi\:D(\xi) \to N$ 
admits a Poincar\'e embedding 
if and only if 
$$
\tilde H^\xi(\tilde c_f) \in  
\{\tilde N/\partial \tilde N, D_2 \tilde P^{f^*\tau_N - \tau_P}\}_\pi
$$
is trivial.
\end{prop}

\section{Immersions of the diagonal}

Suppose $M$ is a closed Poincar\'e space of dimension $d \ge 4$
with fundamental group $\pi$ and stable Spivak tangent bundle $\tau$.
We let $\tilde M\to M$ denote the universal cover. Then
$\pi\times\pi$ is the fundamental group of $M\times M$ and 
$\tilde M \times \tilde M \to M \times M$ is its universal cover.
Consider the $(\pi\times\pi)$-equivariant stable
collapse of the diagonal map $\Delta\: M\to M \times M$:
$$
\tilde c_{\Delta}\in \{\tilde M^+ \smsh \tilde M^+ ,
\hat M^{\tau}\}_{\pi\times \pi} \, .
$$
Here $\hat M^\tau$ is the $(\pi\times \pi)$-equivariant Thom spectrum
defined as follows: consider the regular $(\pi\times \pi)$-covering space
$$
\hat M \to M
$$ 
defined by pulling back the universal cover of $M \times M$.
Then $\hat M^\tau$ is the Thom spectrum of the pullback of $\tau$ to $\hat M$.
Note that $\hat M$ can also be considered as the effect of inducing
$\tilde M$ along the diagonal homomorphism $\pi \to \pi\times \pi$:
$$
\hat M = \tilde M \times_\pi (\pi\times \pi)
$$
and, likewise, $\hat M^\tau = \tilde M^\tau \smsh_\pi (\pi\times \pi)_+$,
where $\tilde M^\tau$ is the Thom spectrum of the pullback of $\tau$
to $\tilde M$.

Now choose such a $(d-1)$-spherical lift $\xi\:S(\xi) \to M$
of $\tau$. Then we have the associated map
$$
\begin{CD}
\Delta_\xi\:D(\xi) \simeq M  @> \Delta >> M \times M\, .
\end{CD}
$$
Applying \ref{better_compression_result}, we immediately get

\begin{prop} \label{xi-embed} The map 
$\Delta_\xi\:D(\xi) \to M \times M$ admits a Poincar\'e embedding
if and only if 
$$
\tilde H^\xi(\tilde c_{\Delta}) \in
\{\tilde M_+ \smsh \tilde M_+ , D_2\hat M^\tau \}_{\pi \times \pi}
$$
is trivial.
\end{prop}

Our next step will be to identify the obstruction group
appearing in \ref{xi-embed}. Let
that $Q_d(\pi) = {\Bbb Z}[\bar \pi]_{\Bbb Z_2}$ 
be the abelian group defined in the introduction.

\begin{thm} \label{isomorphism} 
There is a preferred
isomorphism of abelian groups
$$
\{\tilde M^+ \smsh \tilde M^+ ,
D_2(\hat M^\tau )\}_{\pi\times \pi} \,\, \cong \,\,  Q_d(\pi) \, .
$$
\end{thm}

The proof will make use of the following lemma.

\begin{lem}\label{degree} Let $N$ be a connected 
Poincar\'e duality space
of dimension $n$ equipped with orientation sheaf ${\cal L}$ and fundamental
class $[N]$. Let $E$ be a naive $\pi$-spectrum (or based 
stable $\pi$-space). 
Assume $E$ is  $(n{-}1)$-connected.

Then there is a preferred isomorphism of abelian groups
$$
\{\tilde N^+,E\}_\pi \,\, \cong \,\, H_0(N;{\cal L}\otimes H_n(E))\, .
$$
\end{lem}

\begin{proof} 
Obstruction theory gives an isomorphism
$$
\{\tilde N^+,E\}_\pi   = H^n(N; H_n(E))\, ,
$$
where the right side is the cohomology in degree $n$
of $N$ with coefficients in the $\pi$-module $H_n(E)$.
Poincar\'e duality identifies the right side with 
$$
H_0(N;{\cal L}\otimes H_n(E))\, .
$$
\end{proof}

\begin{proof}[Proof of Theorem \ref{isomorphism}]
The $(\pi\times \pi)$-spectrum $D_2(\hat M^\tau)$ is 
$(2d-1)$-connected,
so by \ref{degree}, we have an isomorphism
$$
\{\tilde M^+ \smsh \tilde M^+ ,
D_2(\hat M^\tau )\}_{\pi\times \pi} \,\, \cong \,\, 
H_0(M\times M;({\cal L}\times {\cal L})\otimes H_{2d}(D_2(\hat M^\tau ))) 
\, ,
$$
where ${\cal L}$ is the orientation bundle for $M$.
Furthermore, we can rewrite the right side as the coinvariants
of the action of ${\Bbb Z_2}$ on
$$
H_0(M\times M; ({\cal L}\times {\cal L})\otimes 
H_{2d}(\hat M^\tau\smsh \hat M^\tau))) 
$$
coming from the self map of
$\hat M^\tau\smsh \hat M^\tau$ which switches factors.

Application of the Thom isomorphism gives
an isomorphism of $(\pi\times \pi)$-modules 
$$
({\cal L}\times {\cal L})\otimes H_{2d}(\hat M^\tau\smsh \hat M^\tau))
\,\, \cong \,\, 
H_0(\hat M) \otimes H_0(\hat M) \, .
$$
 With respect to this
isomorphism, the ${\Bbb Z}_2$-action on $H_0(\hat M) \otimes H_0(\hat M)$
switches factors and multiplies by $(-1)^d$. Consequently,
we are reduced to computing
$$
H_0(M\times M; H_0(\hat M)\otimes H_0(\hat M))
$$
with the given involution. The $(\pi\times \pi)$-module
$H_0(\hat M)$ is ${\Bbb Z}[\pi^\text{ad}]$, where
the latter is the free module on $\pi$ given the
adjoint action of $\pi\times \pi$. 
We therefore have an isomorphism of 
$(\pi\times\pi)$-modules
$$
H_0(\hat M) \otimes H_0(\hat M) \,\, \cong  \,\,
{\Bbb Z}[\pi^{\text{ad}} \times \pi^{\text{ad}}] \, ,
$$
in which the involution on the left corresponds
to the one on the right given by switching factors
and multiplying by $(-1)^d$. Here $\pi^{\text{ad}} \times\pi^{\text{ad}}$
is given the diagonal $\pi\times \pi$-action.

With respect to the above, the abelian group
$$
H_0(M\times M; H_0(\hat M) \otimes H_0(\hat M))
$$
is now identified with the coinvariants of $\pi\times \pi$
acting on ${\Bbb Z}[\pi^{\text{ad}} \times\pi^{\text{ad}}]$.
But this group of coinvariants is just ${\Bbb Z}[\bar \pi]$:
an isomorphism is given by mapping a conjugacy class
$\bar x$ to the pair $(1,x)$, where $x\in\pi$ denotes 
any representative of $\bar x$. With respect to this
isomorphism, the involution on ${\Bbb Z}[\bar \pi]$
is given by $\bar x \mapsto (-1)^d\bar x^{-1}$.
This completes the proof of \ref{isomorphism}.
\end{proof}

\section{Lifts of the tangent fibration}
Let $M$ be a closed, connected
Poincar\'e space of dimension $d \ge 4$ with
Spivak tangent fibration $\tau\: M \to BG$.
Here $BG$ denotes the classifying space for stable spherical
fibrations, and we let $BG_d$ 
denote the classifying space for $(d-1)$-spherical
fibrations. The (inclusion)  map 
$$
i\:BG_d \to BG
$$
is $d$-connected.
So as simplify the exposition, we assume
$i$ is a fibration (by converting it to one).

It follows that
$\tau$ is represented by an $(d-1)$-spherical fibration, i.e.,
there is a section $\xi\: M \to BG_d$ of the fibration $i$ along
$\tau$. The set of {\it lifts} $${\frak L}_\tau$$ 
of $\tau$ is defined to be the
homotopy classes of sections of $i$ along $\tau$. An element
of ${\frak L}_\tau$ is therefore represented by a
$(d-1)$-spherical fibration $\xi$ together with a choice
of stable fiber homotopy equivalence $h\: \xi \simeq_s \tau$.

The fiber of $i$ at the basepoint is identified with
$G/G_d$. It is $(d-1)$-connected, and its $d$-homotopy
group is canonically isomorphic to $Q_d$.

By classical obstruction theory, there is
a free and transitive action
$$
H^d(M;{\cal L}\otimes \pi_d(G/G_d))\times {\frak L}_\tau
\to  {\frak L}_\tau
$$ 
where the cohomology group
$$
H^d(M;{\cal L}\otimes \pi_d(G/G_d)) \, .
$$
encodes the difference obstructions to lifts on the
top dimension of $M$ (recall that ${\cal L}$ is the orientation sheaf).
By Poincar\'e duality, we obtain an isomorphism of the
above with $H_0(M; \pi_d(G/G_d)) = Q_d$.
Summarizing,

\begin{lem} \label{lifts} ${\frak L}_\tau$ is a $Q_d$-torsor,
i.e., there is a free and transitive action 
$$
Q_d \times {\frak L}_\tau \to {\frak L}_\tau\, ,
$$
and therefore, if we choose a basepoint $\xi \in {\frak L}_{\tau}$, 
we obtain an isomorphism of based sets $Q_d \cong {\frak L}_{\tau}$.
\end{lem}

A direct description of the action is as follows:
a result of Wall \cite[2.4]{Wall_PD} shows that
$M$ possesses a top cell decomposition
$$
M \simeq M_0 \cup_\beta D^d\, ,
$$
with $M_0$ a CW complex of dimension $\le d-1$ if $d\ge 4$.
For simplicity, fix
$M = M_0 \cup D^n$. The set of lifts of $\tau$ along
$M_0$ is trivial. Hence, a given
lift $\xi$ is completely described by its value on $D^d$.
Given representative lifts $\xi$ and $\eta$ of $\tau$, we may
assume they coincide on $M_0$. We then form the $(d-1)$-spherical
fibration over $S^d = D_-^d \cup D_+^d$ by taking $\eta$ along the
$D^d_-$ and $\xi$ along $D^d_+$.
This fibration comes equipped with a preferred stable trivialization,
so it defines an element of  $\pi_d(G/G_d) = Q_d$.
Thus pairs of lifts ``differ'' by an element of $Q_d$.

Note the special case when  $\tau$ is trivial gives
a preferred basepoint in ${\frak L}_\tau$, namely, the trivial lift. 
In this case, ${\frak L}_\tau$ is just
$[M,G/G_d]$, and the identification with $Q_d$
is defined by mapping a homotopy class $\gamma\: M \to G/G_d$
to $\gamma_*([M]) \in H_d(G/G_d) = Q_d$.

\begin{rem} The papers of Dupont \cite{Dupont1}, \cite{Dupont2} and Sutherland
\cite{Sutherland} concern explicit detection of the image
of the function
$$
{\frak L}_\tau \to [M,BG_d]\, .
$$
When $d$ is even, the Euler characteristic
$\chi\:[M,BG_d]\to {\Bbb Z}$ detects the image. When $d$ is odd,
the image has either one or two elements, and the invariant
which detects them is  subtler.
\end{rem}

\section{Proof of Theorem \ref{first}}

We will define a function
$$
\phi\: {\frak L}_\tau \to 
\{\tilde M^+\smsh \tilde M^+, D_2\hat M^\tau\}_{\pi\times \pi} \,\,
\overset {\text{\tiny Th.\ } \ref{isomorphism}} \cong \,\, Q_d(\pi) \, .
$$
An element of ${\frak L}_\tau$ is represented by
a spherical fibration and a stable fiber homotopy equivalence
$h\: \xi \simeq_s \tau$. Then $h$ induces an identification
$\hat M^\xi \simeq_s \hat M^\tau$, and the equivariant stable collapse
$\tilde c_\Delta$ becomes a stable map 
$$
\tilde M^+ \smsh \tilde M^+ \stableto \hat M^\xi \, .
$$
Then $\phi(\xi)$ is given by
$$
\tilde H(\tilde c_\Delta) \: \tilde M^+ \smsh \tilde M^+
\stableto D_2 \hat M^\xi \, .
$$

Recall that ${\frak L}_\tau$ comes equipped with
a free and transitive action of $Q_d$. 
Recall that the inclusion
of the trivial group into $\pi$ induces an inclusion
$$
i\:Q_d \subset Q_d(\pi)
$$
of subgroups.

\begin{thm} \label{equivariant} The function $\phi$ is 
one-to-one and has image equal to a coset of $Q_d \subset Q_d(\pi)$.
\end{thm}

\begin{rem} In the category of smooth manifolds, one
can give a short geometric proof of \ref{equivariant} using
Whitney's method \cite[th.\ 3]{Whitney} of locally modifying
the double point set of
a generic immersion by introducing an extra self-intersection. 
This method is
not available in the Poincar\'e category. We will 
use homotopy theoretic methods to
prove \ref{equivariant}.
\end{rem}

We first note how \ref{equivariant} implies
Theorem \ref{first}. By \ref{better_compression_result}, the
diagonal $\Delta\: M \to M\times M$ will admit a Poincar\'e embedding
if and only if the equation $\phi(\xi) = 0$ has a solution. 
By \ref{equivariant}, $\phi({\frak L}_\tau) \subset Q_d(\pi)$
forms a coset of the subgroup $Q_d\subset Q_d(\pi)$. 
Consequently, solutions
to the given equation exist if and only if the coset 
$\phi({\frak L}_\tau)$ contains $0$. Define the 
obstruction
$$
\mu_M \in \tilde Q_d(\pi) \quad (= Q_d(\pi)/Q_d)
$$
to be the equivalence class of the coset $\phi({\frak L}_\tau)$.
Then $\mu_M$ is the complete obstruction to finding a
Poincar\'e embedding of the diagonal. The proof of 
Theorem \ref{first} has now been reduced to establishing \ref{equivariant}.
\medskip

The proof of \ref{equivariant} will require some preparation.
Suppose that $h\in Q_d$
$\xi,\eta \in {\frak L}_\tau$ are elements satisfying
$h\cdot \eta = \xi$. An unraveling of definitions allows use
to think of $h$ as a stable fiber homotopy equivalence
$$
\eta \simeq_s \xi
$$ 
given by composing the pair of 
stable equivalences $\eta \simeq_s \tau \simeq_s \xi$. 
Then $h$ induces a $(\pi\times \pi)$-equivariant stable equivalence 
$$
\hat h\:\hat M^\xi \, \, \stableto \hat M^\eta \, .
$$
The $(\pi\times\pi)$-equivariant stable Hopf invariant $\tilde H$ 
of $\hat h$ is  then an equivariant stable map
$$
\tilde H(\hat h) \: \hat M^\xi \,\, \stableto D_2\hat M^\eta \, ,
$$
which, in view of our identifications, 
we are entitled to regard is an equivariant  stable
map 
$$
\hat M^\tau \,\,  \stableto D_2\hat M^\tau\, .
$$

\begin{lem} \label{action1} With respect to above assumptions,
$$
\phi(\xi) - \phi(\eta)\,\,  =  \,\, 
\tilde H(\hat h)\circ \tilde c_\Delta\, .
$$
\end{lem}

\begin{proof} This is just  a specific instance of
the  composition formula \ref{composition}.
\end{proof}

\begin{lem}\label{action2} The operation
$$
h \mapsto \tilde H(\hat h)\circ \tilde c_\Delta
$$
gives a homomorphism $\psi\:Q_d \to Q_d(\pi)$.
\end{lem}

\begin{proof} This follows again by
the composition formula \ref{composition}. If $g,h \in Q_d$, then
we think of $h$ as stably identifying $\xi$ with $\eta$ and
$g$ as stably identifying $\mu$ with $\xi$. Then 
$$
\tilde H(\hat g\circ \hat h)) = 
\tilde H(\hat g)\circ \hat h + D_2(\hat g) \circ \tilde H(\hat h)
$$
with respect to our identifications, this is the
homomorphism property.
\end{proof}

\begin{lem}\label{action3} The image of
the homomorphism $\psi$ is
contained within the subgroup $Q_d$.
\end{lem}

\begin{proof} A stable fiber homotopy equivalence $h\:\xi \simeq_s \eta$
yields a stable $\pi$-map $\tilde h\: \tilde M^\xi \to \tilde M^\eta$, which
when induced  along the diagonal homomorphism $\pi\to \pi\times \pi$, yields
$\hat h\: \hat M^\xi \to \hat M^\eta$. By \ref{induction}, we have
$$
\tilde H(\hat h) := H_{\pi\times\pi}(\hat h) = i_2 H_\pi(\tilde h) \, ,
$$
which implies that $\psi(h)$ is in the image of
the homomorphism
$$
i_{2*} \: \{\tilde M^+ \smsh \tilde M^+, 
D_2 \tilde M^\tau\smsh_\pi (\pi\times \pi)_+\}_{\pi\times \pi}
\to 
 \{\tilde M^+ \smsh \tilde M^+, 
D_2 \hat M^\tau\}_{\pi\times \pi} 
$$
given by composing with $i_2$.
A calculation of the kind appearing in the
proof of \ref{isomorphism} shows that the domain of
$i_{2*}$ is isomorphic to $Q_d$, 
and that with respect to this
identification, $i_{2*}$ is the inclusion $Q_d \to Q_d(\pi)$.
\end{proof}

Let
$$
\epsilon\: Q_d(\pi) \to Q_d
$$
be the homomorphism induced by mapping 
a conjugacy class to the identity. Note that
$\epsilon$ is a retraction to the inclusion $Q_d \subset
Q_d(\pi)$.

The following lemma implies that $\psi$ is one-to-one.

\begin{lem}\label{action4}
The homomorphism
$$
h \mapsto \epsilon \circ \tilde H(\hat h)\circ \tilde c_\Delta 
$$
is one-to-one.
\end{lem}

\begin{proof} By \ref{induction}, 
the homomorphism has an alternative description
as
$$
h \mapsto H(h^\sharp) \circ c_\Delta \in \{M^+\smsh M^+,D_2 M^\tau\} 
\cong Q_d\, ,
$$
where $h^\sharp \: M^\xi \,\, \stableto M^\eta$ is the 
stable map of Thom spaces associated with $h$, 
$$
H(h^\sharp)\: M^\xi \,\, \stableto D_2 M^\eta  \simeq D_2 M^\tau
$$
its stable Hopf invariant and $c_\Delta\: M^+\smsh M^+ \stableto M^\xi$
the stable collapse. 

By a straightforward argument that we omit,
the homotopy cofiber of the stable map $c_\Delta$ has the weak
homotopy type of a CW spectrum with cells in dimensions $\le 2d-1$.
Since $D_2 M^\tau$ is $(2d-1)$-connected, it follows that
the homomorphism
$$
c^*_\Delta \:\{M^\tau,D_2 M^\tau\} \to \{M^+\smsh M^+,D_2 M^\tau\}
$$
given by composing with $c_\Delta$ is
is surjective. Again, straightforward calculation
shows that the domain and codomain
of this last homomorphism are both 
isomorphic to $Q_d$.  But any surjection of $Q_d$ onto itself is
necessarily an isomorphism. It follows that $c^*_\Delta$ is
an isomorphism.

Consequently, we are reduced to showing that
$$
h \mapsto H(h^\sharp) \in \{M^\tau,D_2 M^\tau\} \cong Q_d
$$
is one-to-one. This is proved below in \ref{last_lemma}.
\end{proof}

\begin{cor} The homomorphism $\psi \: Q_d \to Q_d(\pi)$
is one-to-one with image contained in $Q_d$.
\end{cor}

In view of this corollary, we may to view
$\psi$ as a monomorphism $Q_d \to Q_d$.

\begin{cor} When $d$ is odd, the monomorphism
$\psi\: Q_d \to Q_d$ is an isomorphism.
\end{cor}

\begin{proof} $Q_d$ is finite  when $d$ is odd. 
Since $\phi$ is one-to-one, it must also be onto.
\end{proof}

Now assume that $d$ is even. Then $Q_d$ is infinite
cyclic and we have a monomorphism $2\: Q_d \to {\Bbb Z}$
given by multiplication by $2$.

\begin{lem} \label{action5} Assume $d$ is even. Represent an element
of $Q_d$ by a stable fiber homotopy equivalence
$h\: \eta\simeq_s \xi$
of $(d-1)$-spherical fibrations over $M$. Then up to sign,
$$
2\psi(h) \,\, = \,\, \chi(\xi) - \chi(\eta) \, , 
$$
where $\chi$ is the Euler characteristic.
\end{lem}

\begin{proof} When $d$ is even,
The homomorphism $2\psi$ can alternatively be described as
$$
h \mapsto \text{tr}(H(h^\sharp)) \in 
\{M^\tau,M^\tau\smsh M^\tau\} \cong {\Bbb Z}
$$
where $\text{tr}\:\{M^\tau,D_2 M^\tau\}\to \{M^\tau,M^\tau\smsh M^\tau\}$
is the transfer. By definition of the stable Hopf invariant, 
$\text{tr}(H(h^\sharp))$ coincides with the difference of 
reduced diagonal maps $\Delta_\xi \: M^\xi \to M^\xi \smsh M^\xi$ and 
$\Delta_\eta \: M^\eta\to M^\eta \smsh M^\eta$. But each diagonal 
represents the Euler characteristic of the fibration that it is 
subscripted by.
\end{proof}

\begin{cor} $\psi\: Q_d \to Q_d$ is an isomorphism when $d$ is even.
\end{cor}

\begin{proof} Choose a stable fiber homotopy equivalence
$h\: \xi \simeq_s \eta$ representing a generator of $Q_d$.
We may assume $\xi$ and $\eta$ coincide on $M_0$ (= $M$ with its
top cell $D^d$ removed).
Then $\chi(\xi) - \chi(\eta) = \pm 2$, because one can  
describe the difference obstruction in this case as a map
$$
S^d \to G/G_d
$$
such that the composite $S^d \to G/G_d \to BG_d$
maps the upper hemisphere via $\xi$ and its lower one via $\eta$.
It is well-known that the Euler characteristic of this 
fibration over $S^d$ is $\pm 2$, since the composite
$$
\begin{CD}
Q_d = \pi_d(G/G_d) @>>> \pi_d(BG_d) @>\chi >> {\Bbb Z}
\end{CD}
$$
is multiplication by $2$ when $d$ is even (and trivial when $d$ is odd).

Then \ref{action5} says that $2\psi$ takes value $\pm 2$.
This implies that $\psi$ takes value $\pm 1 \in Q_d$, so
$\psi$ is an isomorphism.
\end{proof}

Recall that if $h\:\xi \simeq_s \eta$ represents an element
of $Q_d$, then $h^\sharp\: M^\xi \simeq_s M^\eta$ denotes
the associated stable equivalence of (unequivariant) Thom spaces.

The next result was used in the proof of \ref{action4}.

\begin{prop} \label{last_lemma} The homomorphism
$$
Q_d \to \{M^\tau,D_2M^\tau\}\, ,
$$
given by 
$
h \mapsto H(h^\sharp)
$
is one-to-one.
\end{prop}

\begin{proof} If $H(h^\sharp)$ is trivial, then
$h^\sharp$ desuspends to an unstable weak equivalence
of Thom spaces 
$$
M^\xi \simeq M^\eta \, .
$$
We need to show that $h\:\xi \simeq_s \eta$ desuspends to
an unstable fiber homotopy equivalence. Let $BF_d$ be
the classifying space of $d$-dimensional spherical fibrations
equipped with section. Since the map
$$
BG_d \to BF_d
$$
given by $\xi \mapsto \xi \oplus \epsilon$ is $(2d-3)$-connected,
and $d \ge 4$, the map $[M,BG_d]\to [M,BF_d]$ is
a bijection. 
It will therefore be enough to show that 
$$
\xi \oplus \epsilon,\eta\oplus \epsilon\: M \to BF_d
$$
are homotopic. 
Let $S^\xi$ and $S^\eta$ denote the total
spaces of these fibrations. Note that $M^\xi$ is just
$S^\xi$ with $M$ collapsed to a point.

Let $B$ be a space equipped with sectioned fibrations $E_1 \to B$
and $E_2 \to B$. Let $[E_1,E_2]_B$ be the fiberwise homotopy classes of
maps $E_1 \to E_2$ which cover the identity of $B$
and preserve the preferred sections, and let
$\{E_1,E_2\}_B$ denote the stable fiberwise
homotopy classes.

Consider the commutative diagram
\begin{equation} \label{EHP}
\xymatrix{
[S^\xi,S^\eta]_M \ar[r]\ar[d] & [M^\xi, M^\eta]\ar[d]\\
\{S^\xi,S^\eta\}_M \ar[r] & \{M^\xi, M^\eta\}\, ,
}
\end{equation}
where the horizontal arrows are given by sending a
fibration map its induced map of Thom spaces, and the vertical ones are
given by stabilization. We also have an identification
$$
[M^\xi, M^\eta] \cong [S^\xi,M^\eta\times M]_M \, ,
$$
so the horizontal arrows of the diagram can also be described
as arising from the $d$-connected map $S^\eta \to M^\eta \times M$ 
of fibrations over $M$ with section by taking homotopy 
classes of maps from $S^\xi$.
The right vertical arrow of \eqref{EHP} sits in a short
exact sequence of based sets whose third term is
$\{M^\xi,D_2 M^\eta\}$, and the map onto this
term is surjective for dimensional reasons (cf.\ \ref{EHPsequence}).
Similarly, the left vertical arrow  
of \eqref{EHP} sits in a short exact sequence whose third
term is  
$$
\{S^\xi,D_2^\bullet S^\eta\}_M \, ,
$$
where $D_2^\bullet S^\eta $ is the quadratic construction applied fiberwise
to $S^\eta$, and the map onto this term is again surjective.

The induced function
\begin{equation} \label{fiberwise}
\{S^\xi,D_2^\bullet S^\eta\}_M \to \{M^\xi,D_2 M^\eta\} \cong Q_d
\end{equation}
arises from a map $D_2^\bullet S^\eta \to D_2^\bullet (M^\eta \times M)$
by applying fiberwise stable maps out of $S^\xi$ 
to the domain and codomain.  The function
\eqref{fiberwise} is therefore a homomorphism.
It is straightforward
to check that the map inducing it is $(2d)$-connected. 
Since $S^\xi$ has the dimension
of a cell complex with cells in dimensions at most $2d$, it follows
that \eqref{fiberwise} is surjective.
\medskip

{\flushleft \it Assertion.} There is an isomorphism 
$$
Q_d \cong \{S^\xi,D_2^\bullet S^\eta\}_M \, .
$$

Here is a proof: first consider the
abelian group
$$
\{S^\xi, S^{2\eta}\}_M
$$
where $S^{2\eta}$ denotes the total space of 
$\eta \oplus \eta \oplus \epsilon$. This has a fiberwise involution
whose homotopy orbits gives  $D_2^\bullet S^\eta$.
Since $\xi$ and $\eta$ are stable fiber homotopy equivalent, we can
identify $\{S^\xi, S^{2\eta}\}_M$ with 
$$
\{S^\eta, S^{2\eta}\}_M \cong \{S^0\times M, S^\eta\}_M \, ,
$$
where this last isomorphism comes from taking fiberwise smash
product with $S^\eta$. The abelian group  
$\{S^0\times M, S^\eta\}_M$ is the fiberwise
homotopy classes of  stable sections
of $S^\eta$. It is also known as the stable cohomotopy of $M^+$
twisted by the spherical fibration $\tau$ (compare \cite[p.\ 5]{Crabb}).
By  Atiyah duality \cite{Atiyah}, it is isomorphic to 
the zeroth stable homotopy group of $M_+$, i.e.,
$$
\pi^{\rm st}_0(M_+)  \cong {\Bbb Z} \, .
$$
Taking into account how the involution acts, we see that 
$\{S^\xi,D_2^\bullet S^\eta\}_M$ is isomorphic to the coinvariants
of a certain involution on ${\Bbb Z}$. It is not difficult to
check that this involution is given by $(-1)^d$. This completes the proof
of the assertion.
\medskip

The proposition is now completed as follows: since 
the map \eqref{fiberwise} is identified with a surjective map
of $Q_d$, it is an isomorphism. Consequently, the diagram
\eqref{EHP} is a cartesian square of abelian groups. The element
$h$ lives in the lower left corner of the square, and by hypothesis,
its image in the lower right corner lifts to the upper right corner.
Consequently, $h$ lifts to unstable fiber homotopy equivalence
$\xi\oplus \epsilon\simeq \xi\oplus \epsilon$ of sectioned fibrations.
The proof of \ref{last_lemma} is now complete.
\end{proof}

\end{document}